%% file: Oscillatory_Sparse_Hilbert_Transform.tex
\documentclass[10pt]{amsart}  
\usepackage[top=1.5in, bottom=1.5in, left=1.25in, right=1.25in]	{geometry}
\usepackage{cite}
\usepackage{mathtools}
\usepackage{parskip}
\usepackage{amsmath,amsthm}
\usepackage[initials]{amsrefs}
\usepackage{amsfonts}
\usepackage{color}
\usepackage{amssymb,bbm}
\usepackage{framed}
\usepackage{multicol}
\usepackage{graphicx}
\usepackage[makeroom]{cancel}
\usepackage{tikz}
\usepackage{hyperref}
\hypersetup{
  colorlinks   = true, 
  urlcolor     = blue, 
  linkcolor    = blue, 
  citecolor   = blue 
}
\usetikzlibrary{arrows}
\input comandos.tex

\allowdisplaybreaks

\theoremstyle{plain} 
\newtheorem{lemma}[equation]{Lemma} 
\newtheorem{proposition}[equation]{Proposition} 
\newtheorem{theorem}[equation]{Theorem} 
\newtheorem{corollary}[equation]{Corollary}

\newtheorem{priorResults}{Theorem}

\theoremstyle{definition}

\theoremstyle{remark}

\numberwithin{equation}{section}

\begin{document}
 
\title[Sparse Discrete Quadratic Hilbert Transform]{Sparse Bounds for Discrete Quadratic Phase
Hilbert Transform}

\author[Kesler]{Robert Kesler}
\author[Mena]{Dar\'io Mena Arias}

\maketitle

\begin{abstract}
Consider the discrete quadratic phase Hilbert Transform acting on $\ell^{2}$ finitely supported
functions 
$$
H^{\alpha} f(n) : = \sum_{m \neq 0} \frac{e^{2 \pi i\alpha m^2} f(n - m)}{m}.
$$
We prove that, uniformly in $\alpha \in \bT$, there is a sparse bound for the bilinear form $\inn{H^{\alpha} f}{g}$. 
The sparse bound implies several mapping properties such as weighted inequalities in an intersection of Muckenhoupt and reverse H\"older classes.
\end{abstract}


\section{Introduction}

Let $e(t) = e^{2 \pi i t}$ and $\alpha \in \bT$.
We consider the operator $H^{\alpha}$ acting on finitely supported functions $f$ on $\bZ$, defined by
$$
H^{\alpha} f(n) : = \sum_{m \neq 0} \frac{e(\alpha m^2) f(n - m)}{m}.
$$
This can be regarded as a discrete oscillatory Hilbert transform with a quadratic phase.  
As such it satisfies a range of $\ell^p$ estimates which are uniform in $\alpha$. In particular, the result below holds. 
Indeed, the work of Arkhipov and Oskolkov \cite{MR0922412} in the case of $p=2$, 
and of Pierce \cite{MR2661174} in the case of $1< p < \infty$, prove much more than the result below. 

\begin{priorResults} 
For $1<p<\infty $, there holds 
$$
\sup _ \alpha \lVert H^{\alpha} : \ell^p \to \ell^p \rVert < \infty . 
$$
\end{priorResults}

In this paper we give a further quantification of the 
uniform boundedness of $H^{\alpha}$, by proving a sparse bound. We set notation for the sparse bound.     Let a discrete interval (or just an interval) be a set of the form $I = \bZ \cap [a,b]$, for $a,b \in
\bR$, and define its length $|I|$ as its cardinality.  For $1 \leq r < \infty$, the $L^r$-average of
a function $f$ on the interval $I$ is defined by
$$
\ave{f}_{I,r} : = \left[  \frac{1}{|I|} \sum_{x \in I} |f(x)|^{r}   \right]^{1/r}.
$$

A collection of intervals $\cS$ is called $\rho$-sparse if for each $S \in \cS$, there is a subset
$E_S$ of $S$ such that (a) $|E_S| > \rho|S|$, and (b) $\|  \sum_{S \in \cS} \cfn{E_S}  \|_{\infty}
\leq \rho^{-1}$.  For a sparse collection $\cS$, a sparse bilinear form $\Lambda$ is defined by
$$
\Lambda_{\cS,r,s}(f,g) : = \sum_{S \in \cS} \ave{f}_{S,r} \ave{g}_{S,r} |S|
$$
When $r= s$, we write $\Lambda_{\cS,r,s} = \Lambda_{\cS,r}$.  The dependence on $\rho$ is not relevant, so it can be omitted.  We also omit sometimes the dependence on the sparse collection  $\mathcal S$ and just write $\Lambda_{r,s}$ or $\Lambda_{r}$.

To simplify some of the arguments, we make use of the following definition:  For an operator $T$
acting on finitely supported functions on $\bZ$, and $1 \leq r, s < \infty$ define its sparse
norm
\begin{equation}\label{SparseNorm}
\| T : \textup{Sparse}(r,s) \| = \| T: (r,s) \|,
\end{equation}
as the infimum over the constants $C > 0$ such that for all finitely supported functions $f,g$ on
$\bZ$ we have
$$
|\inn{Tf}{g}| \leq C \sup \Lambda_{r,s} (f,g).
$$
Here, the supremum is taken over all sparse forms.

With this notation, we can state the main result of this paper as follows,
\begin{theorem}\label{MainResult}
There exists $1 < r < 2$ such that
$$
\sup_{\alpha \in \bT} \|  H^{\alpha} :  (r, r)  \| < \infty.
$$
\end{theorem}

Given the useful structure of the sparse forms, we can derive a variety of mapping properties.  For instance, we obtain the following immediate result
\begin{corollary}
There exists $1 < r < 2$ such that for all weights $w$ that satisfy $w, w^{-1} \in A_2 \cap RH_r$
we have
$$
\|  H^{\alpha} : \ell^2(w) \mapsto \ell^2(w)  \| \lesssim 1.
$$
\end{corollary}
The weights above are in the intersection of the the standard Muckenhoupt class $A_2$ and some  Reverse H\"older class $RH_r$.  Here and through all the paper, the notation $A \lesssim B$ means that there is a constant $C$ such that $A \leq C  B$;  the dependence of the constant will be indicated when necessary.

The domination by sparse operators has been an active topic initiated by Lerner \cites{MR3085756} in his simple proof of the $A_2$ conjecture, by providing sparse control over the norm of a Calder\'on-Zygmund operator.  This was improved to a pointwise estimate in \cites{MR3521084} and  following a stopping time argument in \cites{2015arXiv150105818L}.  The latter approach has been used in different contexts \cites{2015arXiv151000973B, 2015arXiv151005789H, MR3484688}.  The sparse bilinear form approach that we use here, has proven to be successful where the pointwise approach is not convenient or to avoid the use of maximal truncations, for example, the bilinear Hilbert transform \cites{2016arXiv160305317C}, Bochner-Riesz multipliers \cites{2016arXiv160506401B}  and oscillatory singular integrals \cites{2016arXiv160908701K, 2016arXiv160906364L}.

The study of oscillatory singular integrals is motivated by the work of Stein, who in \cites{MR1364908} proves the boundedness on $L^p$, for $1 < p < \infty$, of the following operator,
\begin{equation}\label{SteinOscillatoryQuadratic}
\sup_{\alpha \in \bR} \left|  \int_{\bR} f(x-y) \frac{e(\alpha y^2)}{y} \, dy  \right|.
\end{equation}

In the setting of discrete norm inequalities it is important to mention the remarkable work of Bourgain on ergodic theorems regarding polynomial averages \cites{MR937581,MR937582}.  More recent results include the work of Krause \cites{2014arXiv14021803K}  which have been extended in different directions by 
Mirek, Stein and Trojan \cites{2016arXiv151207518M, 2016arXiv151207523M}.  A first result in which similar discrete operator can be controlled by sparse forms can be found in \cites{2016arXiv161208881C}, and in the case of random discrete operators in \cites{2016arXiv160906364L,2016arXiv160908701K}, where the sparse bound follows from simpler arguments.

Our main result, and the proof, is a model case for a wider range of results in the discrete setting.   Some of the many possible extensions to the main result of this paper are as follows. 
\begin{enumerate}
\item Extend the result to a general polynomial and kernel.  That is, given a polynomial $P$ and a Calder\'on-Zygmund Kernel $K$, find sparse bounds for the operator
$$
T_P f(n) = \sum_{m \neq 0} e(P(m))K(m)f(n-m),
$$
that only depend on the degree of $P$ and the kernel.  More ambitious claims suggest themselves, such as obtaining sparse bounds for discrete Radon transforms, even in the quasi-translation invariant setting.   See \cites{MR2661174,MR3053569}.  
\item Sparse version of Krause and Lacey's result \cite{2015arXiv151206918K}, that is, find sparse bounds for the following restricted maximal operator, for $A$ satisfying a certain Minkowski dimension condition,
$$
\sup_{\alpha \in A} \left|H^{\alpha} f(n) \right|  =  \sup_{\alpha \in A} \left| \sum_{m \neq 0} \frac{e(\alpha m^2) f(n - m)}{m} \right|.
$$
\item Sparse control over the maximal truncations of the  operators above. This would entail extra difficulties.  

\end{enumerate}

The paper is organized as follows: In \S 2, we provide some preliminary results regarding sparse forms and specific operators bounded by them, that are key to our proof.  In \S 3, following techniques from the Hardy-Littlewood circle method, we give a decomposition for the Fourier multiplier of the operator into major and minor arc components, and obtain some estimates for the different parts.  We prove the sparse bounds for the minor and major arcs in \S 4 and \S 5 respectively to conclude the main theorem.  Of particular interest is the method to bound the major arcs, as it depends upon the sparse bound in Theorem \ref{SparseOscillatory}.


\section{Preliminaries}

One useful fact about sparse operators is that, in some sense, they admit an universal domination. 
A version of the following lemma can 
be found in \cite{2016arXiv161001531L} and has a similar proof.
\begin{lemma}\label{UniversalDomination}
Given finitely supported functions $f,g$ and $1 \leq r, s < \infty$, there is a sparse form
$\Lambda^{*}_{r,s}$ and a constant $C> 0$ such that for any other sparse operator $\Lambda_{r,s}$
we have
$$
\Lambda_{r,s}(f,g) \leq C \Lambda_{r,s}^{*}(f,g).
$$
\end{lemma}

The Hardy-Littlewood maximal function is defined by
\begin{equation}\label{HLMaximalFunction}
\cM_{\textup{\tiny HL}} f  (n) : = \sup_{N \geq 0} \frac{1}{2N + 1} \sum_{j = -N}^N |f(n-j)|, \qquad n \in \bZ.
\end{equation}
A well known result is the following. 
\begin{priorResults}\label{MHLSparse}
The Hardy-Littlewood maximal function satisfies $(1,1)$ sparse bounds.  That is,
$$
\|  \cM_{\textup{\tiny HL}} : \textup{Sparse}(1,1) \| \lesssim 1.
$$
\end{priorResults}

If $\cH$ is a Hilbert space, we extend the definition of sparse forms to vector valued functions $f$, by setting $\ave{f}_I = |I|^{-1}\sum_{x \in I} \| f(x) \|_{\cH}$.   It is then straightforward to extend some sparse domination results to Hilbert space valued functions.  One of this results, in the continuous setting of oscillatory singular integrals, is the following theorem, that is going to be an important part of our proof.  
\begin{priorResults}\label{SparseOscillatory}\cites{2016arXiv160906364L}
Let $K$ be a Calder\'on-Zygmund kernel and $P$ a polynomial of degree $d$ on $\bR^n$.  Define the operator 
$$
T_P f(x) = \int_{\bR} e(P(y)) K(y) f(x - y) \, dy.
$$
For each $1 <  r < 2$ and compactly supported, Hilbert space valued functions $f, g$, there is a constant $C = C(K,d,n,r)$ and a bilinear sparse form $\Lambda_r$ such that
$$
\| T_P f  : \textup{Sparse}(r,r) \| \leq C.
$$

\end{priorResults}
Recall that a Calder\'on-Zygmung kernel $K : \bR \bs \{ 0 \} \en \cC$ satisfies 
$$
\sup_{y \neq 0}|y K(y)| + \left|y^2 \tfrac{d}{dy} K(y)\right| < \infty,
$$
and the corresponding convolution operator is $L^2(\bR^n)$-bounded.  In particular, we are going to apply this result with the Hilbert Transform kernel $K(y) = 1/y$.  It is important to note that the previous estimate depends on the polynomial only through its degree.

In the subsequent sections, $\eps > 0$ will denote a small fixed constant.  We use the standard
notations for the Fourier transform and its inverse:
\begin{align*}
\hat{f}(\beta)   & = \cF f (\beta)  = \sum_{n \in \bZ} f(n) e(-\beta n), \\
\check{g}(n)   & = \cF^{-1} g (n)  = \int_{\bT} g(\beta) e(\beta n) \, d\beta. \\
\end{align*}


\section{Decomposition of the multiplier}

The Fourier multiplier associated to the transformation $H^{\alpha}$ is
\begin{equation}\label{multiplier}
M^{\alpha}(\beta) : = \sum_{m \neq 0} \frac{e(\alpha m^2 - \beta m ) }{m}.
\end{equation}
The  goal of this section is to describe a decomposition of the multiplier $M^\alpha$ into terms, with uniform control in the variable $\alpha$.
Let $\{ \psi_j \}_{j \geq 0}$ be a dyadic resolution of the function $\frac{1}{t}$, with $\psi_j(t)
= 2^{-j}\psi( 2^{-j} t)$, and $\psi$ is a odd smooth function satisfying $\psi(t) \leq \cfn{ [1/4,
1] }(|t|)$.  Then, for $|t| \geq 1$, we have $\frac{1}{t} = \sum_{j \geq 0} \psi_j(t)$, and in the
support of $\psi_j$, we have $2^{j-2} \leq |t| \leq 2^j$. Using this,
we can decompose the multiplier as a sum of terms of the form
$$
M_j^{\alpha} (\beta) : = \sum_{m \neq 0} \frac{e(\alpha m^2 - \beta m ) }{m} \psi_j(m).
$$
That way, we can write $M^{\alpha} = \sum_{j} M_j^{\alpha}$.

For fixed $s \in \bN$, define
$$
\cR_s : = \set{ \left( \tfrac{A}{Q}, \tfrac{B}{Q}  \right) \in \bT^2  }{  A,B,Q \in \bZ,\ 
(A,Q)=(B,Q)=1 ,\ 2^{s-1} \leq Q  \leq 2^{s}}.
$$
Then, the rationals in the torus, can be written as $\bigcup_{s \in \bN} \cR_s$.  Given
$(\frac{A}{Q}, \frac{B}{Q}) \in \cR_s$, and $j \geq s/ \eps$, define the $j$-th major arc at
$(\frac{A}{Q}, \frac{B}{Q})$ by

\begin{equation}\label{majorbox}
 \mathfrak{M}_j(A/Q, B/Q) := \left\{ (\alpha, \beta) \in \mathbb{T}^2: |\alpha - A/Q| \leq
2^{(\epsilon-2) j}, |\beta - B/Q| \leq 2^{(\epsilon-1)j} \right\}.
\end{equation}
Collect the major arcs 
\begin{equation}\label{majorboxcollection}
 \mathfrak{M}_j  : = \bigcup_{\substack{ (A,Q)=(B,Q)=1 \\ 0< Q \leq 2^{6 \eps j}}}
\mathfrak{M}_j(A/Q, B/Q).
\end{equation}
As proven in \cite{2015arXiv151206918K}, the union above is over disjoint sets for $\eps$ small enough.  For each $j$, we define the minor arcs to be the complement of this union of the major arcs.

Let $\chi$ be a smooth even bump function, such that $\cfn{[-1/10, 1/10]} \leq \chi \leq \cfn{[-1/5,
1/5]}$.  For $s,j \in \bN$, set $\chi_s(t): = \chi(10^{s}t)$, and define the multiplier
\begin{equation}\label{Ljs}
L_{j,s}^{\alpha}(\beta) : = \sum_{\left(\frac{A}{Q}, \frac{B}{Q}\right) \in \cR_s} S(A/Q, B/Q)
U_j(\alpha -A/Q,\beta-B/Q ) \chi_s(\alpha -A/Q)\chi_s(\beta -B/Q).
\end{equation}
Here, $U_j$ is a continuous analogue of the multiplier $M_j$,
\begin{equation}\label{ContinuousAnalogue}
U_j (x, y) : = \int_{\bR} e(xt^2 - yt) \psi_j(t) \, dt,
\end{equation}
and $S$ is the complete Gauss sum
\begin{equation}\label{CompleteGaussSum}
S(A/Q, B/Q) : = \frac{1}{Q} \sum_{r = 0}^{Q-1} e(  A/Q \cdot r^2 - B/Q \cdot r).
\end{equation}

Consider also the following definitions.
\begin{align}
L_j^{\alpha}(\beta) & : = \sum_{s \leq  \eps j} L_{j,s}^{\alpha}(\beta), \quad j \geq 1, \\
L^{\alpha,s}(\beta) & : = \sum_{j \geq  s/ \eps} L_{j,s}^{\alpha}(\beta), \quad s \geq 1, \label{Ls}
\\
L^{\alpha}(\beta) & : = \sum_{j = 1}^{\infty} L_j^{\alpha}(\beta) = \sum_{s = 1}^{\infty}
L^{\alpha,s}(\beta), \\
E_j^{\alpha}(\beta) & : = M_j^{\alpha}(\beta) - L_j^{\alpha}(\beta), \quad j \geq 1, \label{Ej} \\
E^{\alpha}(\beta) & : = \sum_{j = 1}^{\infty} E_j^{\alpha}(\beta).
\end{align}

The proof of the following lemmas can be found in \cite{ 2015arXiv151206918K }.  The first one says
that on the major arcs, $M_j$ is well approximated by its continuous analogue.  
\begin{lemma}
For $1 \leq s \leq \eps j$, $( A/Q, B/Q ) \in \cR_s$, and $(\alpha, \beta) \in \mathfrak{M}_j(A/Q,
B/Q)$, we have the approximation
$$
M_j^{\alpha}(\beta) = S(A/Q, B/Q) U_j(\alpha -A/Q,\beta-B/Q ) + O(2^{(3\eps -1)j}).
$$
\end{lemma}

In the minor arcs we have the following estimates.
\begin{lemma}
There exists $\delta = \delta(\eps)$ such that uniformly in $j \geq 1$,
$$
|M_j^{\alpha}(\beta)| + |L_j^{\alpha}(\beta))| \lesssim 2^{-\delta j}, \quad (\alpha,\beta) \not \in
\mathfrak{M}_j(A/Q, B/Q)
$$
\end{lemma}

Using these results, we obtain the following bounds.
\begin{theorem}
There is a choice of $\delta > 0$ such that, uniformly in $\alpha \in \bT$
\begin{align}
|S(A/Q, B/Q)| & \lesssim 2^{-\delta s}, \quad (A/Q, B/Q) \in \cR_s, \quad s \geq 1, \label{GaussSumEstimate} \\
\| E_j^{\alpha}(\beta) \|_{\infty} &\lesssim 2^{-\delta j}, \quad j \geq 1, \label{EjLInfEstimate}\\
\left\| \pdd{}{\beta}{\beta}E_j^{\alpha}(\beta) \right\|_{\infty} &\lesssim 2^{2j}, \quad j \geq 1.
\label{EjDerivEstimate}
\end{align}
\end{theorem}

The first estimate can be found in several places in the literature (see, for example, \cite{Hua}).
Given that by construction $M_j^{\alpha}(\beta) = L_j^{\alpha}(\beta) + E_j^{\alpha}(\beta)$, the
second estimate is a consequence of the previous two lemmas.  The derivative estimate comes from
straightforward computations.

We prove the main Theorem by showing that there is a choice of $1 < r < 2$ and $\eta > 0$ such that for $j,s \geq 1$ the following estimates hold, uniformly in $\alpha \in \bT$
\begin{align}
\|  T_{\check{E}_j^{\alpha}}  : (r,r) \| & \lesssim 2^{-\eta j} \qquad \mbox{(Minor
arcs estimate)} \label{MinorArcsEstimate} \\
\|  T_{\check{L}^{\alpha,s}}  : (r,r) \| & \lesssim 2^{-\eta s} \qquad \mbox{(Major
arcs estimate)} \label{MajorArcsEstimate}
\end{align}
Since our operator can be written as $H^{\alpha}  = \sum_{j \geq 1} T_{\check{E}_j^{\alpha}}  + \sum_{s \geq 1} T_{\check{L}^{\alpha,s}}$, from the triangle inequality for the sparse norm it follows that
$$
\| H^\alpha : (r,r) \| \leq  \sum_{j \geq 1} \|  T_{\check{E}_j^{\alpha}}  : (r,r) \| + \sum_{s \geq 1} \|  T_{\check{L}^{\alpha,s}}  : (r,r) \| \leq \sum_{j \geq 1} 2^{-\eta j} + \sum_{j \geq 1} 2^{-\eta s} < \infty.
$$
Since these estimates are independent of $\alpha$, the main theorem follows.


\section{Minor Arcs estimate} 

Consider the multiplier $E_j^{\alpha}$, defined in (\ref{Ej}).  The $L^{\infty}$ estimate
(\ref{EjLInfEstimate}) and the derivative estimate (\ref{EjDerivEstimate})
imply that
\begin{equation}\label{boundonKernelEj}
|\cF^{-1}E_j^{\alpha}(m)| \lesssim \min\left\{ 2^{-\eps j}, \frac{2^{2j}}{1+m^2} \right\}.
\end{equation}
These bounds are independent of $\alpha$, since the derivative estimates are.  Write
$\cF^{-1}E_j^{\alpha} = \check{E}_{j,1}^{\alpha} +
\check{E}_{j,2}^{\alpha}$, where $\check{E}_{j,1}^{\alpha}(m) = \cF^{-1}E_j^{\alpha}(m)
\cfn{[-2^{3j}, 2^{3j}]}(m)$.  We first estimate for $\check{E}_{j,2}$,  for this, consider the
Hardy-Littlewood maximal function defined in (\ref{HLMaximalFunction}), we have
\begin{align*}
|T_{\check{E}_{j,2}^{\alpha}} f (x)| & = |\check{E}_{j,2}^{\alpha} * f (x)| \leq \sum_{y \in \mathbb{Z}} |K_2(y)
f(x-y)| \lesssim 2^{2j} \sum_{|y| \geq 2^{3j}}  \frac{|f(x-y)|}{1+|y|^{2}} \\
& = 2^{2j} \sum_{|k| \geq 3j} \sum_{2^{k} \leq |y| < 2^{k+1}} \frac{|f(x-y)|}{1+|y|^{2}}
\lesssim 2^{2j} \sum_{|k| \geq 2^{3j}} 2^{-k} \cM_{\textup{\tiny HL}} f(x) = 2^{-j} \cM_{\textup{\tiny HL}} f(x).
\end{align*}
Once again, this estimate is independent of $\alpha$.  Using Theorem \ref{MHLSparse} we obtain the
result for $\check{E}_{j,2}^{\alpha}$.

For $\check{E}_{j,1}^{\alpha}$, we need to use the following result (Proposition 2.4 in
\cite{2016arXiv161208881C}).
\begin{proposition}\label{p:elementary}
Let $ T _{K} f (x) = \sum_{n} K (n) f (x-n)$ be convolution with kernel $ K$.  Assuming that $ K$ is
finitely supported on the interval $ [-N,N]$ we have  the  inequalities 
$$
\lVert T_K  : (r,s)\rVert  \lesssim  N ^{1/r+1/s-1} \lVert T _{K}  :  \ell ^r
\mapsto  \ell ^{s'}\rVert , \qquad 1\leq r, s < \infty . 
$$
\end{proposition}

To proof the sparse bound for $T_{\check{E}_{j,1}^{\alpha}}$, we use the proposition with $N =
2^{3j}$ and $r = s$, that is
$$ 
\|  T_{\check{E}_{j,1}^{\alpha}}  : (r,r) \| \lesssim 2^{3j(\frac{2}{r} - 1)}  \lVert
T_{\check{E}_{j,1}^{\alpha}}  :  \ell ^r \mapsto  \ell ^{r'} \rVert, \qquad 1 \leq r < \infty.
$$

We just need to find an $r$ such that the operator norm has a summable decay in $j$.  It is easy to
check for the cases $r= 1$ and $r = 2$. For $r = 1$, we have by Young's
inequality and (\ref{boundonKernelEj})
$$
\left\|  T_{\check{E}_{j,1}^{\alpha}} f  \right\|_{\infty} \lesssim \|  \check{E}_{j,1}^{\alpha}
\|_{\infty} \| f \|_{1} \lesssim 2^{-\delta j} \| f \|_{1}
$$
And for $r = 2$, we have by the $L^{\infty}$ estimate of the multiplier
\eqref{EjLInfEstimate}, and Plancherel,  
$$
\left\lVert T_{\check{E}_{j,1}^{\alpha}}  :  \ell ^2 \mapsto  \ell ^{2} \right\rVert \lesssim
2^{-\delta j}.
$$
We can now interpolate and choose $1 < r < 2$ such that $10(2/r-1)  < \delta/2$ to get the desired
decay.  Combining this with the estimate over the norm of $T_{\check{E}_{j,2}^{\alpha}}$ the proof
of \eqref{MinorArcsEstimate} is complete.


\section{Major Arcs estimate}

We proceed now to prove the more complicated estimate \eqref{MajorArcsEstimate}.  Recall the
definition of $U_j$, given by \eqref{ContinuousAnalogue}.  For $s\geq 0$, define $U^s$ to
be
$$
U^s(x,y) = \sum_{j \geq s/\eps} U_j(x,y).
$$
Then, we can write the multiplier $L^{\alpha, s}$ defined in (\ref{Ls}) as
$$
L^{\alpha, s}(\beta) = \sum_{\left(\frac{A}{Q}, \frac{B}{Q}\right) \in \cR_s} S(A/Q, B/Q) U^s(\alpha
-A/Q,\beta-B/Q ) \chi_s(\alpha -A/Q)\chi_s(\beta -B/Q).
$$
Given that the support of $\chi_s$ is contained in $\left[-2\cdot 10^{-s-1},2\cdot
10^{-s-1}\right]$, for fixed $\alpha \in \bT$, there is at most one rational $\alpha_s = A/Q$ with
$(A,Q)=1,
2^{s-1} \leq Q \leq 2^{s}$ for which $\chi_s( \alpha -A/Q )$ is non zero.  To simplify the notation,
we make use of the following definition
$$
\cR_s^{\alpha} = \set{ B/Q \in \bT }{ (A/Q, B/Q) \in \cR_s, A/Q = \alpha_s }.
$$
It is important to say that the subsequent analysis only depends upon the cardinality of $\mathcal{R}_s ^\alpha$, which is at most $2^{2s}$, and not the value of $\alpha$.  We can rewrite $L^{\alpha, s}$ as
$$
L^{\alpha, s}(\beta) = \sum_{\frac{B}{Q} \in \cR_s^{\alpha}} S(\alpha_s, B/Q) U^s(\alpha
-\alpha_s,\beta-B/Q ) \chi_s(\alpha -\alpha_s)\chi_s(\beta -B/Q).
$$

As in \cite{2016arXiv161208881C}, we will make use of a sparse bound for Hilbert space valued
singular integrals.   For this, define for fixed $\alpha \in \bT$ the finite dimensional Hilbert space $\cH_s^{\alpha} = \ell^2 ( \cR^{\alpha}_s )$.  Given $f \in \ell^2$,  if
$\textup{Mod}_{h} f(x)= e(hx)f(x)$ represents the standard modulation by $h$, set $f_{s,h} : =
\cF^{-1}(\chi_s^{1/2}) * \textup{Mod}_{-h} f$.   Define the $\cH_s^{\alpha}$-valued function
$f_s^{\alpha}$ by
$$
f_s^{\alpha}  : = \set{f_{s,B/Q}}{ B/Q \in \cR^{\alpha}_s }.
$$
Note that the Fourier transforms $\hat{f}_{s,B/Q}(\beta) = \chi_s^{1/2} (\beta) \hat{f}(\beta +
B/Q)$ have disjoint supports, so by Bessel's Theorem $ \| f_s^{\alpha} \|_{\ell^2(\cH_s)} \leq \| f
\|_{\ell^2}$.  We have the following simplifications,
\begin{align*}
\inn{T_{\check{L}^{\alpha, s}}f}{g} & = \sum_{\tfrac{B}{Q} \in \cR_s^{\alpha}} \sum_{j \geq s/\eps} S(\alpha_s,\tfrac{B}{Q}) \inn{U_j(\alpha-\alpha_s,\cdot-\tfrac{B}{Q})\chi_s(\alpha-\alpha_s)\chi_s(\cdot-\tfrac{B}{Q})\hat{f}(\cdot)}{\hat{g}(\cdot)} \\
& = \sum_{\tfrac{B}{Q} \in \cR_s^{\alpha}} \sum_{j \geq s/\eps} S(\alpha_s,\tfrac{B}{Q}) \inn{U_j(\alpha-\alpha_s,\cdot)\chi_s(\alpha-\alpha_s)\chi_s^{1/2}(\cdot)\hat{f}(\cdot + \tfrac{B}{Q})}{\chi_s^{1/2}(\cdot)\hat{g}(\cdot + \tfrac{B}{Q})} \\
& = \chi_s(\alpha-\alpha_s) \sum_{\tfrac{B}{Q} \in \cR_s^{\alpha}} \sum_{j \geq s/\eps} S(\alpha_s,\tfrac{B}{Q}) \inn{U_j(\alpha - \alpha_s, \cdot)\hat{f}_{s,B/Q}}{\hat{g}_{s,B/Q}}  \\
& =  \chi_s(\alpha-\alpha_s) \sum_{\tfrac{B}{Q} \in \cR_s^{\alpha}} S(\alpha_s,\tfrac{B}{Q}) \inn{ T_{\check{U}^s} f_{s,B/Q}}{g_{s,B/Q}}.
\end{align*}

For $B/Q \in \cR_s^{\alpha}$, take $\lambda_{B/Q}$ with unit norm and such that $\lambda_{B/Q} \inn{ T_{\check{U}^s} f_{s,B/Q}}{f_{s,B/Q}} \geq 0$, and set $\tilde{f}_s^{\alpha} = \set{\lambda_{B/Q} f_{s,B/Q}}{ B/Q \in \cR_s^{\alpha} }$, then, using the Gauss sum estimate (\ref{GaussSumEstimate}) and summing over $j\geq s/\eps$ we have
$$
|\inn{T_{\check{L}^{\alpha, s}}f}{g}| \lesssim 2^{-\delta s} \inn{T_{\check{U}^s} \tilde{f}_s^{\alpha} }{g_s^{\alpha}}.
$$
Since $\| f_s^{\alpha} \|_{\cH_s^{\alpha}} = \| \tilde{f}_s^{\alpha} \|_{\cH_s^{\alpha}}$, then we can replace $\tilde{f}_s^{\alpha}$ by ${f}_s^{\alpha}$ in the inner product.   The next step is to find a sparse form $\Lambda_1$ (on Hilbert space valued functions) that dominates the last inner product.  For this, first we write
\begin{align*}
U^s(\alpha - \alpha_s,\beta) & = \sum_{j \geq s/\eps} \int_{\bR} e((\alpha - \alpha_s)t^2) e( -\beta t ) \psi_j(t) \, dt = \int_{\bR} e((\alpha - \alpha_s)t^2 - \beta t) \sum_{j \geq s/\eps} \psi_j(t) \, dt
\end{align*}
The integrand above is supported on $|t| \geq 2^{\lfloor s/\eps \rfloor -2}$ and by explicit computation $\sum_{j \geq s/\eps} \psi_j(t)$ coincides with $\frac{1}{t}$ for $|t| \geq 2^{\lfloor s/\eps \rfloor}$.  Therefore, this kernel corresponds to a Calder\'on-Zygmund kernel, and we can apply Theorem \ref{SparseOscillatory}. 
As a consequence,  for any $1 < r_1 < 2$  there is a sparse bilinear form $\Lambda_{r_1}$ such that 
\begin{equation}\label{VSpaceSparse}
|\inn{T_{\check{U}^s}f_s}{g_s}| \lesssim \Lambda_{r_1}(f_s^{\alpha},g_s^{\alpha}).
\end{equation}
 The implied constant above  does not depend on $\alpha$.

To end the proof, we need the following result
\begin{lemma}
Let $1 \leq r_1 <2$ and $\delta >0$. Let $\Lambda_{r_1}$ be a sparse form over a collection of intervals all of which have length larger than $10^s$. Then there exists $r$ satisfying $r_1< r < 2$ such that for all $f,g$ there is a sparse form $\Lambda_{r}$ for which
$$
\Lambda_{r_1}( f_s^{\alpha}, g_s^{\alpha} ) \lesssim 2^{\delta s/4} \Lambda_{r}(f,g).
$$
\end{lemma}
The proof of this lemma is a slight modification of the proof of the $r_1=1$ result given at the end of \cite{2016arXiv161208881C} (the value of $\alpha$ doesn't affect the proof). Ensuring all the sparse intervals in $\Lambda_{r_1}$ have length at least $10^s$ is achieved by taking $\epsilon>0$ small enough.  Combining the estimates, and letting $\eta = 3\delta / 4$, we have
$$
|\inn{T_{\check{L}^{\alpha, s}}f}{g}| \lesssim 2^{-\eta s} \Lambda_r(f,g).
$$
Which proves the major arcs estimate, and therefore, the main theorem.

\bibliographystyle{amsplain}

\end{document}

%% file: comandos.tex
\newcommand{\bc}{\begin{center}}
\newcommand{\ec}{\end{center}}
\newcommand{\beqn}{\begin{align*}}
\newcommand{\eeqn}{\end{align*}}
\newcommand{\benu}{\begin{enumerate}}
\newcommand{\eenu}{\end{enumerate}}
\newcommand{\bit}{\begin{itemize}}
\newcommand{\eit}{\end{itemize}}

\newcommand{\cfn}[1]{\ensuremath \raisebox{0pt}{$\mathbbm{1}$}_{#1}}

\newcommand{\en}{\rightarrow}

\newcommand{\bN }{\mathbb{ N }}

\newcommand{\bR }{\mathbb{ R }}
 
\newcommand{\bT }{\mathbb{ T }}

\newcommand{\bZ }{\mathbb{ Z }}

\newcommand{\cC }{\mathcal{ C }}

\newcommand{\cF }{\mathcal{ F }}
 
\newcommand{\cH }{\mathcal{ H }}

\newcommand{\cM }{\mathcal{ M }}

\newcommand{\cR }{\mathcal{ R }}
\newcommand{\cS }{\mathcal{ S }}

\newcommand{\eps}{\varepsilon}

\newcommand{\bs}{\backslash}

\newcommand{\pdd}[3]{\ifx#2#3\frac{\partial^2 #1}{\partial #2^2}\else \frac{\partial^2 #1}{\partial #2\,\partial #3}\fi}

\newcommand{\ave}[1]{\left\langle #1 \right\rangle}
\newcommand{\inn}[2]{\left\langle #1 , #2 \right\rangle}

\newcommand{\set}[2]{\left\{\,#1 : #2 \, \right\}}

